\documentstyle[amstex,12pt]{amsart}
\renewcommand{\baselinestretch}{1.3}

\title{Multi-graded extended Rees algebras of $\m$-primary ideals}
\author{Clare D'Cruz} 
\date{}

\makeatletter
\newcommand{\single}{\let\CS=\@currsize\renewcommand{\baselinestretch}{1}\tiny\CS}
\newcommand{\singles}{\let\CS=\currsize\renewcommand{\baselinestretch}{1.3}\tiny\CS}
\newcommand{\oneanda}{\let\CS=\currsize\renewcommand{\baselinestretch}{1.2}\tiny\CS}
\newcommand{\doubles}{\let\CS=\currsize\renewcommand{\baselinestretch}{1.5}\tiny\CS}
\newcommand{\tree}{\let\CS=\currsize\renewcommand{\baselinestretch}{1.5}\tiny\CS}
\newcommand{\four}{\let\CS=\currsize\renewcommand{\baselinestretch}{2}\tiny\CS}

\newcommand{\ncom}{\newcommand} 
\ncom{\bq}{\begin{equation}}
\ncom{\eq}{\end{equation}}
\ncom{\beqn}{\begin{eqnarray*}}
\ncom{\eeqn}{\end{eqnarray*}}
\ncom{\beq}{\begin{eqnarray}}      
\ncom{\eeq}{\end{eqnarray}}
\ncom{\been}{\begin{enumerate}}
\ncom{\eeen}{\end{enumerate}}
\ncom{\nno}{\nonumber}
\ncom{\rar}{\rightarrow}
\ncom{\lrar}{\longrightarrow}
\ncom{\Rar}{\Rightarrow}
\ncom{\noin}{\noindent}

\newtheorem{thm}{Theorem}[section]
\newtheorem{lemma}[thm]{Lemma}
\newtheorem{cor}[thm]{Corollary}
\newtheorem{pro}[thm]{Proposition}
\newtheorem{example}[thm]{Example}
\newtheorem{remark}[thm]{Remark}
\newtheorem{blank}[thm]{}

\ncom{\bt}{\begin{thm}}
\ncom{\et}{\end{thm}}
\ncom{\bl}{\begin{lemma}}
\ncom{\el}{\end{lemma}}
\ncom{\bco}{\begin{cor}}
\ncom{\eco}{\end{cor}}
\ncom{\bp}{\begin{pro}}
\ncom{\ep}{\end{pro}}
\ncom{\bex}{\begin{example}}
\ncom{\eex}{\end{example}}
\ncom{\brm}{\begin{remark}}
\ncom{\erm}{\end{remark}}
\ncom{\bb}{\begin{blank}}
\ncom{\eb}{\end{blank}}
\ncom{\bc}{\begin{center}}
\ncom{\ec}{\end{center}}

\ncom{\comx}{{\mathbb{C}}}
\ncom{\ze}{{\mathbb{Z}}}
\ncom{\re}{{\mathbb{R}}}
\ncom{\Q}{{\mathbb{Q}}}
\ncom{\N}{{\mathbb{N}}_0}

\ncom{\CM}{Cohen-Macaulay }
\ncom{\f}{\frac}
\ncom{\al}{\alpha}
\ncom{\be}{\beta}
\ncom{\ga}{\gamma}
\ncom{\bib}{\bibitem}
\ncom{\sms}{\setminus}
\ncom{\seq}{\subseteq}

\ncom{\olin}{\overline}
\ncom{\bip}{\bigoplus}
\ncom{\sta}{\stackrel}

\def\b{{\cal B}}
\def\r{{\cal R}}
\def\n{{\cal N}}
\def\i{{\bf I}}
\def\j{{\cal J}}
\def\m{{\frak m}}

\begin{document} 
\maketitle
\section{Introduction}

In this paper we consider multi-graded extended Rees algebras of zero
dimensional ideals which are Cohen-Macaulay (CM) with minimal
multiplicity. We show that the minimal multiplicity property can occur
only for  the ordinary extended Rees algebra and the bigraded extended
Rees algebra. For the bigraded extended Rees algebra we find
necessary conditions   for it to be CM with minimal
multiplicity. 
  We also produce bigraded 
Rees algebras which are Cohen-Macaulay with minimal multiplicity.  

 A considerable amount was known for the ordinary extended Rees
 algebra. Among the many we quote  (\cite{katz-ver}, \cite{verma2},
 \cite{verma3}, \cite{verma6}). There was nothing known concerning the
 minimal multiplicity of the multi-graded extended Rees algebra.  One
 of the crucial results needed was the formula of multiplicity of a
 maximal homogeneous ideal. This formula was obtained in  the author
 in \cite{clare2}.

Throughout this paper $(R, \m)$ will denote a Noetherian local ring of
positive dimension. Without loss of generality we will assume that
$R/ \m$ is infinite. It is well-known that for any CM local ring $(R,
\m)$, $e(\m) \geq \mu(\m) - \dim~R +1$, where $e(\m)$ denotes the
multiplicity  of $\m$ and $\mu(\m)$ is the minimal number of
generators of $\m$. A CM local ring is said to have {\em minimal
  multiplicity} if equality holds.

 Let $I_{1}, \ldots, I_{g}$ be ideals of positive height in $(R, \m)$
 and let $t_{1}, \ldots, t_{g}$ be indeterminates. The {\em
 multi-graded extended Rees algebra} of $R$ with respect to the ideals
 $I_{1}, \ldots, I_{g}$ is the graded ring $\b(\i) := \oplus_{r_j \in
 \ze, 1 \leq j \leq g} (I_{1} t_{1})^{r_{1}} \cdots (I_{g}
 t_{g})^{r_{g}}$. Here $I_j^{r_{j}} = R$, if $r_j \leq 0$ for all $j =
 1, \ldots,g$.  Let $\n$ denote the maximal homogeneous ideal of
 $\b(\i)$.  The {\em multi-Rees algebra} is the graded ring
 $\bip_{r_{i} \geq 0} (I_{1} t_{1})^{r_{1}} \cdots (I_{g}
 t_{g})^{r_{g}}$ and will be denoted by ${\r}({\i})$.

In the past decade several researchers have investigated the
multi-Rees algebra. Since the  multi-Rees algebra is a subring of the
multi-graded extended Rees algebra, it is natural to expect them to
have similar ring-theoretic properties.  However, there was no
progress concerning  the multi-graded extended Rees algebra. 

Hence we will briefly state some of the earlier known results on the
Rees algebra and the extended Rees algebra.  It is well-known that if
$I$ is an ideal of positive height in a CM local ring $R$ and if
$\r(I)$ is CM, then the {\em associated graded ring} $G(I) :=
\oplus_{r \geq 0}I^r/I^{r+1}$ is also CM
\cite[Proposition~1.1]{hu}. It is easy to see that $G(I)$ is CM if
and only if the extended Rees ring $\b(I)$ is.  In 1989, Verma showed
that if $R$ is a CM local ring of dimension two with minimal
multiplicity, then for all positive integers $r$, $\r(\m^r)$ and
$\b(\m^r)$ are CM with minimal multiplicity
\cite[Theorem~3.3,~4.3]{verma2}.  In the same year he showed that if
$I$ is 
a parameter ideal in a CM ring of dimension at  least two and if $\ell
(I + \m^2 / \m^2) \geq \dim~R-1$, then $\r(I)$ and $\b(I)$ are CM with
minimal multiplicity \cite[Theorem~3.1,~3.2]{verma3}.  In 1991 he
proved the following: Let $(R, \m)$ be a regular local ring of
dimension two. Let $I$ be a contracted $\m$-primary ideal with
reduction number one. Then $\r(I)$ and $\b(I)$ are CM with minimal
multiplicity \cite[Theorem~3.1,~4.3]{verma6}.

In \cite{hhrt} Herrmann et. al. remarked that if  $I$ is an ideal of
positive height and if $I_1 = \cdots = I_g = I$,  then 
 the multi-Rees algebra ${\r}({\i})$ behaves like the ordinary
 extended  Rees
algebra ${\r}(I)$. In this paper they studied the CM property of the multi-Rees
algebra. Minimal multiplicity of the multi-Rees algebra has   been
studied in \cite{verma7}, \cite{hhrt} and    \cite{clare3}.


The following results  which were obtained in the authors thesis
played an important role in obtaining our results:

\been
\item
A  relation between the number of generators of an $\m$-primary ideal
in a CM local ring and a certain mixed multiplicity (Theorem~\ref{nog}).

\item
The bounds on the mixed multiplicities of ideals in $(R, \m)$
(see Lemma~\ref{kv2}, Lemma~\ref{isw}). 

\item
The bounds on   
$\ell(I_1 + I_2 + \m^2/ \m^2)$ when  $\b(\i)_\n$ 
is CM with minimal multiplicity, where $I_1$ and $I_2$ are
$\m$-primary ideals in a CM ring $(R, \m)$ (see Lemma~\ref{con.m.sq},
Lemma~\ref{length}).  
\eeen

\brm
{\em The above mentioned results also  give a simple and unified
  proof  for the known results for the ordinary  extended Rees
  algebra. We do not mention these results here.   But we   answer a
  question of Verma concerning the ordinary extended Rees   algebra
  (see   \cite[pg 3015]{verma2} and Example~\ref{ex2}). This   gives
  an infinite class of  examples of ordinary extended Rees   algebras
  which are Cohen-Macaulay with minimal multiplicity even the 
  original ring does not have minimal multiplicity. It was not
  possible to easily see or  construct this example with the methods
  used in Verma's paper concerning the ordinary extended Rees algebra.  
 }
\erm

 We now summarise the main results in  this paper. 
In Section two we prove that for a CM local ring minimal multiplicity
can occur only for the ordinary extended Rees algebra  and the
bigraded extended Rees algebra (i.e. when $g = 1,2$).  
In Section three we obtain necessary
conditions for the bigraded extended Rees algebra to be CM with
minimal multiplicity. In Section four we consider  bigraded
extended  Rees algebras which are CM with minimal multiplicity.  
We end the paper with an example.

\noin
{\bf Acknowledgements:} The author is very grateful to J. Verma for his
valuable suggestions. The author  also wishes to thank the National
Board for Higher Mathematics  for financial
support and the Indian Institute of Technology, Bombay,  where
research on this paper was carried out. 


\section{Preliminaries}
\label{prelim} 
\noin
\bb
{\em An ideal $J \seq I$ is a {\em
reduction}\index{reduction} of $I$ if there exists a positive integer
$r$ such that $J I^r = I^{r+1}$ \cite{north-rees}. The ideal $J$ is
called a {\em minimal reduction}\index{minimal reduction} of $I$ if
$J$ is minimal  with respect to inclusion among all reductions of
$I$. If $R/\m$ is infinite, then any  minimal reduction of $I$ is
generated by $a(I)$ elements, where $a(I) = \dim \bip_{_{n \geq 0}}
(I^n/\m I^n)$ is called {\em the analytic spread of
$I$}\index{analytic spread}. For an ideal $I$ in $R$, $ht\; I \leq
a(I) \leq \dim \; R$ \cite{rees2}. If $J$ is a reduction of $I$, then
the reduction number of $I$ with respect to $J$ is defined to be  
$$
r_J (I) = \min\{ n \geq 0 \; | \; JI^n = I^{n+1} \}.
$$
The {\em reduction number}\index{reduction number} of $I$ is defined to be 
$$
r(I) = \min \{ r_J (I) \; | \; J \;\; 
               \mbox{is a minimal reduction of } I \}.
$$}
\eb

\bb
{\em Let $I_{1}, \ldots, I_{g}$ be $\m$-primary 
ideals in a local ring $R$ of dimension $d$. For 
$r_{1}, \ldots, r_{g}$ large, 
$\ell_R (I_{1}^{r_{1}}   I_{2}^{r_{2}} \cdots I_{g}^{r_{g}}/
         I_{1}^{r_{1}+1} I_{2}^{r_{2}} \cdots I_{g}^{r_{g}})$ 
is a polynomial of degree $d-1$ in $r_1, \ldots, r_g$ and 
can be written in the form
$$
 \sum_{q_{1} + \cdots + q_{g} = d-1}
 e(I_{1}^{[q_{1}+1]} | I_{2}^{[q_{2}]} |\cdots | I_{g}^{[q_{g}]})
  \f{r_{1}^{q_{1}}}{q_{1}!} \cdots \f{r_{g}^{q_{g}}}{ q_{g}!}    
+ \mbox{lower degree terms};
$$
where
$e(I_{1}^{[q_{1}+1]} | I_{2}^{[q_{2}]} |\cdots | I_{g}^{[q_{g}]})$ 
are positive integers and they are called the {\em mixed
multiplicities} of the set of ideals $\{I_{1}, \ldots, I_{g}\}$  \cite{teis}.
For  $g=2$, we will use the notation 
$$
e_{q}(I_{1} | I_{2}) := e(I_{1}^{[d-q]} | I_{2}^{[q]})
\hspace{1in} 0 \leq q \leq d-1.
$$}
\eb

\bb
{\em Rees obtained an interpretation
of mixed multiplicities in terms of joint reductions \cite{rees5}. Let
$I_{1}, \ldots, I_{g}$ be $\m$-primary ideals.   A set of elements $x_{1},
\ldots, x_{g}$ is a called a {\em joint reduction} of the set of
ideals $\{I_{1}, \ldots, I_{g}\}$
if $x_j \in I_j$ for $j=1, \ldots, g$ and if 
$\sum_{j=1}^{g} x_{j} I_{1} \cdots \widehat{I_{j}} \cdots {I_g}$
is a reduction of $I_{1} \cdots I_{g}$. 
Rees proved that if $R/\m$ is infinite, $g = \dim~R$ and $I_{1}, \ldots
I_{g}$ are $\m$-primary ideals,  then joint reductions exist \cite{rees5}. 
It follows that if $I$ and $J$ are $\m$-primary ideals in a local ring
$(R, \m)$ then $e_0(I|J) = e(I)$ \cite{r1}. We end
this section by stating an  important result of Rees.}
\eb

\bl
\label{rees}
{\bf (Rees' Lemma)} {\em \cite[Lemma 1.2]{rees5}} Suppose
$(R,\m)$ is a local ring with infinite residue field. Let $\{I_{1},
\ldots , I_{g}\}$ be a set of  ideals of $R$ and let  $\cal P$ be
a finite collection of prime ideals of $R$ not containing any of $I_{1},
\ldots , I_{g}$. Then for each $i = 1, \ldots, g$, there exists an
element $x_i \in I_i, \; x_i$ not contained in any prime ideal
of $\cal P$ and an integer $s_i$ such that for $r_i \geq s_i$
and for all positive integers $r_{1}, , \ldots, \widehat{r_i},
\ldots, r_{g}$; 
\beqn
 x_iR \cap I_{1}^{r_{1}} \cdots I_{g}^{r_{g}}
            = x_i I_{1}^{r_{1}} \cdots I_i^{r_i-1} \cdots I_{g}^{r_{g}}.
\eeqn
\el

\section{The case $g \geq 3$}

The main result in this section is:

\bt 
\label{geq3}
Let $(R, \m)$ be a CM local ring of dimension $d$. Let $I_{1}, \ldots,
I_{g}$  be $\m$-primary ideals in $R$. If $\b(\i)_{\n}$ is
CM with minimal multiplicity, then $g \leq 2$. 
\et

By  a result of Valla,   $\dim~\b(\i) = \dim~R +g$
\cite{valla}. Notice that  $e(\n \b(\i)_\n) = e(\n)$ and 
$\mu({\n}\b(\i)_{\n}) = \mu(\n)$.  Hence, if $\b(\i)_{\n}$ is CM, then
it has minimal multiplicity if and only if 
 $e({\n}) = \mu(\n) - (\dim~R + g) + 1$.

We  state an interesting
inequality which relates the  number of generators of an ideal with a
certain mixed multiplicity. 

\bt
\label{nog} 
Let $(R, \m)$ be  a CM local ring of positive dimension $d$. Let $I$ be an 
$\m$-primary ideal  of $R$. Then 
$$
\mu(I) \leq e_{d-1}(\m|I) +  d-1.
$$
\et
\pf We induct on $d$. The case $d=1$ has been proved by J.~Sally
\cite[pg 49]{sally2}.  If $d>1$, then by Lemma~\ref{rees} there exists a
non-zero-divisor $y \in I$ and a positive integer $s_0$ so that
for $s \geq s_0$ and  $r>0$,  
$$
yR \cap \m^r I^s = y \m^r I^{s-1}.
$$
Let ``-'' denote the image in $\olin{R} = R/yR$.
By induction hypothesis we have
\beqn
\mu(I) \leq \mu(\olin{I}) + 1
       \leq e_{d-2}(\olin{\m} \, | \, \olin{I})+ d-1 
       =  e_{d-1}(\m | I) + d-1.
\eeqn
{\qed}

An  upper bound on the number of
generators of the maximal homogeneous ideal of the multi-graded extended
Rees algebra can be estimated by Theorem~\ref{nog}. 

\brm
\label{remone}
{\em Let $I_{1}, \ldots,  I_{g}$ be $\m$-primary ideals in a CM local
ring $(R, \m)$.  Put $L = I_{1} + \cdots + 
I_{g} + \m^2$. Comparing the graded components of $\n$ and $\n^2$ we get   
\beq 
\label{rem}
\nno
\mu(\n)
&=&     \sum_{j=1}^{g} \ell (R/ \m) + \ell (\m/L)
      + \sum_{j=1}^{g} \mu(I_{j}) \\ \nno
&=&    g + \mu(\m)  + \sum_{j=1}^{g} \mu(I_{j}) 
      - \ell \left(L/\m^2 \right)\\
&\leq&   e(\m) + \sum_{j=1}^{g} e_{d-1}(\m|I_{j}) + d(g+1) -1 
     - \ell \left(L/\m^2 \right)~~~~~~~~~~~~~~~~~~~~~
                                       \mbox{[by Theorem~\ref{nog}]}.
\eeq} 
\erm

The multiplicity of $\n$ can be expressed in terms of mixed
multiplicities of ideals in $R$. Hence,  any bound on mixed
multiplicities of ideals in $R$ will  give a bound on the
multiplicity of $R$.

\bt
\cite[Theorem~1.2]{clare2}
\label{mult}
Let $I_{1}, \ldots, I_{g}$ be $\m$-primary ideals in $(R, \m)$. Put 
$L = I_{1} +  \cdots + I_{g} + \m^2$. Then
\beqn
   e(\n) 
= \f{1}{2^d} \left[ 
  \sum_{n=0}^{g} 
  \sum_{q=0}^{d-1} 2^{d-1-q}
  \sum_{\sta{\scriptstyle  q_{1}+ \cdots + q_n = d-1-q}
       {1 \leq i_{1} <  \cdots < i_n \leq g}}
   e({L^{[q+1]} | I_{i_{1}}^{[q_{1}]} | 
     \cdots | I_{i_n}^{[q_{n}]}}) \right]  .
\eeqn
\et

\bl
\label{kv2}
Let $(R,\m)$ be a  local ring and  $I_{1}, \ldots, I_{g}$ be 
$\m$-primary ideals of $R$. Then for all nonnegative integers
$q_{1}, \ldots,  q_{g}$ satisfying  $q_{1} + \cdots + q_{g} = d-1$, 
$$
       e(I_{1}^{[q_{1}+1]}| I_{2}^{[q_{2}]} |\cdots | I_{g}^{[q_{g}]})
 \geq   e( I_{1} +  \cdots + I_{g}).
$$
\el
\pf Since $I_{1}, \ldots, I_{g}$ are $\m$-primary, by
\cite[Theorem~2.4]{rees5}, there exists a joint reduction $x_{1},
\ldots, x_d$ of $q_{1}+1$ copies of $I_{1}$, $q_{2}$ copies of $I_{2}$, \ldots,
$q_{g}$ copies of $I_{g}$, such that $e(x_{1}, \ldots, x_d) =
e(I_{1}^{[q_{1}+ 1]}| I_{2}^{[q_{2}]}\ldots | I_{g}^{[q_{g}]})$.  
Since $(x_{1}, \ldots, x_d) \seq I_{1} + \cdots + I_{g}$, 
$e(x_{1}, \ldots, x_d) \geq e(I_{1} +  \cdots + I_{g})$. 
\qed

\bl
\label{isw}
{\em [cf. Sw, Lemma~2.8]}
Let $(R,\m)$ be a local ring. Let $I_{1}, \ldots, I_d$ be
$\m$-primary ideals in $R$. Let $x_i \in I_i$ for $i=1, \ldots,
d$ be such that $(x_{1}, \ldots, x_d)$ is $\m$-primary. Then 
$$
e ( I_{1} | \cdots |I_d) \leq e(x_{1}, \ldots, x_d). 
$$
If $(R,\m)$ is quasi-unmixed and equality holds, then $x_{1},
\ldots, x_{d}$ is a joint reduction of the set of ideals $\{I_{1},
\ldots, I_{d}\}$.  
\el

\noin
{\bf Proof of Theorem~\ref{geq3}:} 
Put  $L =  I_{1} + \cdots + I_{g} + \m^2$. 
Since $\b(\i)_{\n}$ is CM with minimal multiplicity, $e(\n) = \mu(\n)
- \dim~\b (\i) + 1$. From Remark~\ref{remone} it follows that 
\beq
\label{aaa}
     \mu(\n) - \dim~\b (\i) + 1 
&\leq& e(\m) + \sum_{j=1}^{g}e_{d-1}(\m|I_{j}) + g(d-1) -\ell (L/\m^2)
       \nno \\
&\leq& e(\m) + \sum_{j=1}^{g}e_{d-1}(\m|I_{j}) + g(d-1).
\eeq

Let  $d=1$. Then
\beq
\label{mu}
&&   e(\n)  
\leq (g+1) e(\m)  
      \hspace{1.8in}\mbox{[from (\ref{aaa})]}\\
\label{one}
\mbox{and}  
&&  e(\n) 
=   2^{g-1} e(L)  \geq 2^{g-1} e(\m)  
     \mbox{\hspace{1in} [from Theorem~\ref{mult}]}.
\eeq
Clearly, $2^{g-1} > g+1$ for $g > 3$. 
 If $g=3$, then equality holds in (\ref{mu}) and (\ref{one}). This
implies  
 that $e(L) = e(\m)$ and $L = \m^2$ which is not possible.
Hence $g \leq 2$. 

Let  $d \geq 2$ and $g \geq 3$.  It is enough to show that 
\beq
\label{star}
e(\n) > e(\m) + \sum_{j=1}^{g}e_{d-1}(\m|I_{j}) + g(d-1).
\eeq
Since $I_{1}, \ldots, I_{g}$ are
$\m$-primary ideals, the mixed multiplicities which appear in the
formula of $e(\n)$ (see Theorem~\ref{mult})
are positive integers. Moreover,  
 $e_{d-1}(L|I_{j}) \geq
e_{d-1}(\m|I_{j})$ for all $j = 1, \ldots, g$ (Lemma~\ref{isw}). 
In the multiplicity formula for $e(\n)$, if we replace  $e_{d-1}(L|I_{j})$ by 
$e_{d-1}(\m|I_{j})$
($1 \leq j \leq g$) and  the remaining terms by 1 we get 
\beq
\label{mult11}
e(\n) \nno
&\geq& \f{1}{2^d} 
       \left[ 1+  \sum_{n=1}^{g} {g \choose n} 
       \left[
       \sum_{q=0}^{d-1} 2^{q} {q + n-1 \choose n-1} 
   -                         2^{d-1}  n \right] \right] 
   +    2^{g-2}\sum_{j=1}^{g} e_{d-1}(\m|I_{j}) \\
&=&   \f{1}{2^d} 
      \left[ 1+  \sum_{n=1}^{g} {g \choose n}
      \sum_{q=0}^{d-1} 2^{q} {q + n-1 \choose n-1} \right]
  -    2^{g-2} g 
  +    2^{g-2}\sum_{j=1}^{g} e_{d-1}(\m|I_{j}).
\eeq
 Clearly
\beq
\label{mult12} \nno
&&       2^{g-2} \sum_{j=1}^{g} e_{d-1}(\m|I_{j})  -  2^{g-2} g \\ \nno
&\geq& \sum_{j=1}^{g} e_{d-1}(\m|I_{j}) 
   +     e(\m) 
   +    (2^{g-2} -1)g - 1
   -    2^{g-2} g
       \hspace{1cm} \mbox{[by Lemma~\ref{kv2}]}\\
&=&    \sum_{j=1}^{g} e_{d-1}(\m|I_{j}) 
   +    e(\m) - g - 1.
\eeq

We will show by induction on $d$ that
\beq
\label{mult13}
1 + \sum_{n=1}^{g} {g \choose n} 
       \sum_{q=0}^{d-1} 2^q {q + n-1 \choose n-1} 
>   2^d (gd +1). 
\eeq


If $d=2$,  then it is easy to see that the left-hand side of 
(\ref{mult13}) is $2^{g}(g+1)$ and  the right-hand side is 
$4(2g+1)$. If  $d \geq 3$, then 
\beqn
&&   1 + \sum_{n=1}^{g} {g \choose n} 
         \sum_{q=0}^{d-1} 2^q {q + n-1 \choose n-1} \\
&=&  1 + \sum_{n=1}^{g} {g \choose n} 
         \sum_{q=0}^{d-2} 2^q {q + n-1 \choose n-1}
       +      2^{d-1} 
         \sum_{n=1}^{g} {d-1 + n-1 \choose n-1} 
                        {g \choose n} \\
&>&           2^{d-1} [g(d-1) +1] + 
       2^{d-1} \sum_{n=1}^{g} {g \choose n}{d-1  + n-1 \choose n-1}
               \hspace{1cm}           
               \mbox{[by induction hypothesis]} \\
&>&    2^{d-1} \left[ g(d-1) + 
       1 + {g \choose 1} + 
           {g \choose 2} d\right] \\
&>&    2^d (gd +1).
\eeqn
Comparing (\ref{mult11}), (\ref{mult12}) and (\ref{mult13}) we get the
inequality in (\ref{star}). This completes the proof of the theorem. 
\qed

\section{The Case $g=2$}

In this section we  obtain necessary conditions for the   bigraded
extended Rees algebra to be  CM with minimal multiplicity.  

Let $I_1$ and $I_2$ be ideals of positive height in $R$. Put $L = I_1
+ I_2 + \m^2$. Recall that if  $g = 2$ then 
\beq
\label{geq2}
\label{two1} \nno
&&     e(\n) \\ \nno
&=&   \f{1}{2^d} \left[ e(L) 
    + \sum_{q=0}^{d-1} 
      \sum_{j=1}^{2}2^q e_{q}(L|I_{j})  
    + \sum_{q=0}^{d-1} 2^{d-1-q}
      \sum_{q_{1} + q_{2} =d-1-q} 
       e(L^{[q+1]} |I_{1}^{[q_{1}]} | I_{2}^{[q_{2}]})  \right] \\ \nno 
&\geq& \f{e(L)}{2^d} \left[ 
        1 + \sum_{n=1}^{2} {2 \choose n} \left[ 
        \sum_{q=0}^{d-1} 2^q {q + n-1 \choose n-1} 
     - 2^{d-1} n \right] \right] 
     + e_{d-1} (L|I_{1}) + e_{d-1}(L|I_{2}) \\
&=&   (d-1)e(L) + e_{d-1} (L|I_{1}) + e_{d-1}(L|I_{2}).
\eeq

Putting  $g=2$ in Remark~\ref{remone}, we get
\beq
\label{two3} \nno
&&      \mu(\n) - \dim~\b(\i) + 1 \\ \nno 
&=&     \mu(\m) + \mu(I_{1}) + \mu(I_{2}) 
      - \ell(L/\m^2) - (d -1) \\
&\leq&   e(\m) +  e_{d-1} (\m|I_{1}) 
               +  e_{d-1}(\m|I_{2}) 
               +  2(d-1) 
      - \ell(L/\m^2).
\eeq

\bl
\label{con.m.sq}
Let $I_{1}$ and $I_2$ be $\m$-primary ideals in a CM local ring  $R$
of positive dimension $d$.
If $\b(\i)_{\n}$ is CM with minimal multiplicity, then 
$\ell(I_{1}+I_{2} + \m^2 / \m^2) >0$. 
\el
\pf Suppose  not. Then $I_{1}+I_{2} \seq \m^2$.
It is easy to see that $e(\m^n) = n^d e(\m)$
and $e_{q}(\m^r|I_{j}) = r^{d-q} e_{q}(\m|I_{j})$ for all $q = 1,
\ldots, d-1$ for $j = 1,2$.  Hence from (\ref{two1})    
\beq
\label{two4}  \nno
          e(\n)
&\geq&    (d-1) e(\m^2) +   e_{d-1}(\m^2|I_{1})  
                        +  e_{d-1}(\m^2|I_{2}) \\ 
& =  &    2^d(d-1)  e(\m) 
       +  2  e_{d-1}(\m|I_{1})  +  2e_{d-1}(\m|I_{2}).
\eeq
Since $I_1 + I_2 \seq \m^2$, from (\ref{two3}) we get 
\beq
\label{two5}
\mu(\n) - \dim~\b(\i) + 1 
\leq e(\m) + e_{d-1}(\m|I_{1}) + e_{d-1}(\m|I_{2}) +  2(d-1).
\eeq
Our assumption on  $\b(\i)_\n$ implies that  
$ e(\n) = \mu(\n) - \dim~\b(\i) + 1$.
Hence from (\ref{two4}) and (\ref{two5}) we get 
\beqn
    [ 2^{d} (d-1) -1] e(\m) + e_{d-1}(\m|I_{1}) + e_{d-1}(\m|I_{2}) 
\leq 2(d-1).
\eeqn
Observe that 
$ 2^{d} (d-1) + 1 > 2(d-1)$
for all $d \geq 1$.
This leads to a contradiction. 
\qed

\bl
\label{length}
Let $(R, \m)$ be a CM local ring of dimension $d \geq 2$. Let $I_1$
and $I_2$ be $\m$-primary ideals in $R$. If $\b(\i)_{\n}$ is CM with
minimal multiplicity, then  $\ell(I_{1} + I_{2} + \m^2/\m^2) \leq
d$. If $d \geq 3$, then equality holds.  
\el
\pf Put $L = I_{1} + I_{2} + \m^2$. Since $\b(\i)_{\n}$ is CM with minimal
multiplicity, $e(\n) = \mu(\n) - \dim~\b(\i) + 1$. From
Lemma~\ref{isw},  $e_q(L|I_j) \geq e_q(\m |I_j)$ for $j =
1,2$. Hence from 
 (\ref{geq2}) and (\ref{two3}) we get 
\beqn
(d-2) e(\m) \leq  2(d-1) - \ell(L/\m^2).
\eeqn
Since $e(\m) \geq 1$,   $\ell(L/\m^2) \leq  d$.
Let $d \geq 3$. Assume that 
$\ell(L/ \m^2) \leq d-1$. Then $e(L) \geq e(\m) + 1$. 
Once again from     (\ref{geq2}) and (\ref{two3}) we get
\beq
\label{two6} \nno
&&       2d-3 
      +  e(\m) 
      +  e_{d-1}(\m|I_{1}) 
      +  e_{d-1}(\m|I_{2})\\ \nno
&\leq&  (d-1) [e(\m) + 1] 
      +  e_{d-1}(\m|I_{1}) 
      +  e_{d-1}(\m|I_{2})\\ \nno
&\leq&  (d-1) e(L) 
      +  e_{d-1}(L|I_{1}) 
      +  e_{d-1}(L|I_{2}) \\ 
&\leq&   e(\m) + 2(d-1)
      +  e_{d-1}(\m|I_{1}) 
      +  e_{d-1}(\m|I_{2}) 
      -  \ell(L/\m^2).
\eeq
This gives  $\ell(L / \m^2) \leq 1$. Lemma~\ref{con.m.sq}
implies that $\ell(L / \m^2)=1$. 
Put   $\ell(L/\m^2) = 1$ in (\ref{two6}). Then equality holds in   
(\ref{two6}) and hence
$e(\m) = 1$ and   $e(L) = 2$. Thus
\beqn
     \mu(\m) 
=    \ell \left( \f{\m}{\m^2} \right)
=    \ell \left( \f{R}{L} \right) + \ell \left( \f{L}{\m^2} \right)
   - \ell \left( \f{R}{\m} \right)
=    \ell \left( \f{R}{L} \right)
\leq e(L) 
=    2.
\eeqn 
This leads to a contradiction.
Hence $\ell(L/\m^2) = d$.
\qed


\bl
\label{ii}
Let $(R,\m)$ be a local ring of positive dimension $d$. Let
$I_{1}$ and $I_{2}$ be ideals of positive height in $R$. If 
${\cal B}(I_{1},I_{2})_{\n}$ is CM with
minimal multiplicity, then $r(I_{1}) \leq 1$ and $r(I_{2}) \leq 1$. 
\el
\pf Let
$J_{i}$  be  a minimal reduction of $I_{i}$, ($i = 1,2$). 
Then
$\j = (t_{1}^{-1}, t_{2}^{-1}, \m, J_{1}  t_{1}, J_{2}  t_{2})$ is a 
reduction of $\n$. 
Since $\b(I_{1}, I_{2})_{\n}$ is CM with minimal multiplicity, 
 $J \n = \n^2$ \cite[Theorem~1]{sally1}. Comparing the graded components 
of $\j \n$ and $\n^2$ we get $J_{1}I_{1} + \m I_{1}^2 = I_{1}^2$ and
$J_{2}I_{2} + \m I_{2}^2 = I_{2}^2$. By Nakayama's lemma, 
$J_{1} I_{1} = I_{1}^2$ and $J_{2} I_{2} = I_{2}^2$.
\qed

We are now ready to prove the main results of this section. 

\bt
\label{last}
Let $(R, \m)$ be a CM local ring of dimension $d \geq 3$. Put $L = I_1
+ I_2 + \m^2$.   
Suppose  $\b(\i)_{\n}$ is CM with minimal
multiplicity. Then  
\been
\item
$R$ is a regular local ring;

\item
For $j = 1,2$:
\been
\item 
$\mu(I_j) = e_{d-1}(\m | I_{j})  + d-1$;

\item
$e_{q}(L|I_{j}) = 1$  for all 
$q = 0, \ldots d-2$; 

\item
$r(I_{j}) \leq 1$. 

\eeen
\eeen
\et
\pf Put $L = I_1 + I_2 + \m^2$.  Recall that 
\beq
\label{two7} \nno
&&        e(\n) \\ \nno
&=&    \f{1}{2^d} \left[ e(L) 
     + \sum_{q=0}^{d-1} 2^q \left[ e_{q}(L|I_{1}) + e_{q}(L|I_{2}) \right] 
     + \sum_{q=0}^{d-1} 2^{d-1-q}
       \sum_{q_{1} + q_{2} = d-1-q} 
             e(L^{[q + 1]} | I_{1}^{[q_{1}]} | I_{2}^{[q_{2}]}) \right]\\ 
&\geq& (d-1) e(\m) 
+            e_{d-1}(\m|I_{1}) + e_{d-1}(\m|I_{2}) 
\eeq
and 
\beq
\label{two8} 
\nno
&&
        \mu(\n) - \dim~\b(\i) + 1  \\ \nno
&\leq& e(\m) + e_{d-1}(\m|I_{1}) 
             + e_{d-1}(\m|I_{2}) + 2(d-1) - \ell(L/ \m^2)\\
&=&    e(\m) + e_{d-1}(\m|I_{1}) 
             + e_{d-1}(\m|I_{2}) + d-2 
             \hspace{1in}\mbox{[by Lemma~\ref{length}]}.
\eeq 
Since $\b(\i)_{\n}$ is CM with minimal multiplicity, $e(\n) = \mu(\n) -
\dim~\b(\i) + 1$. Hence from (\ref{two7}) and (\ref{two8}) we get
$(d-2) e(\m) \leq d-2$.
This implies that $e(\m) =1$. Hence  equality holds in (\ref{two7}) and
(\ref{two8}). As a consequence for $j = 1,2$ and  
$q = 0, \ldots, d-2$ we have that 
$\mu(I_j) = e_{d-1}(\m | I_{j})  + d-1$
and $e_{q}(L|I_{j}) = 1$. 
By Lemma~\ref{ii},  $r(I_{j}) \leq 1$ ($j = 1,2$). 
\qed

\bt
\label{last1}
Let $(R, \m)$ be a CM local ring of dimension $d = 2$. Assume that
$\ell (I_{1} + I_{2} + \m^2/\m^2) = 2$.  
Suppose  $\b(\i)_{\n}$ is CM with minimal
multiplicity. Then  
\been
\item
$R$ has minimal multiplicity;

\item
For $j = 1,2$:
\been
\item
 $\mu(I_j) = e_{1}(\m | I_{j})  + 1$;

\item
$r(I_{j}) \leq 1$.  
\eeen
\eeen
\et
\pf The proof of the theorem is similar to the proof of  Theorem~\ref{last}.
\qed

\section{Special cases and Examples}

We recall a result on minimal multiplicity. 

\brm
\label{sally1}
{\em \cite{verma3}, (2.3)
Let $(R,\m)$ be a $d$-dimensional local ring. If $R$ satisfies the
equation of minimal multiplicity, then $R$ is CM if and only if
$J \m= \m^2$ for some minimal reduction $J$ of $\m$.}
\erm

\bt
\label{5.2}
Let $(R, \m)$ be a CM local ring of positive dimension $d$. 
 Let $r$ be a positive
integer.  
\been
\item
If $d =1$, then $\b(\m,\m^{r})_{\n}$ is CM with minimal
multiplicity if and only if $R$ has minimal multiplicity.  
\item
If $d=2$,
 then $\b(\m,\m^{r})_{\n}$ is CM with minimal  multiplicity if and
 only if $R$ is a regular local ring. 
\item
If $d \geq 3$, then
 $\b(\m,\m^{r})_{\n}$ is CM with minimal multiplicity if and only if
 $R$ is a regular local ring and $r=1$.  
\eeen
\et
\pf The necessary part   can be easily verified  for $d
=1$. If  $d = 2$,   it follows from
Theorem~\ref{last1}(2a). Let $d \geq 3$. 
Since $\b(\m,\m^{r})_{\n}$ is CM with minimal multiplicity,  by 
Theorem~\ref{last}, $R$ is a  regular local ring and  
$\mu(\m^r) = e_{d-1}(\m|\m^r) + d-1 = r^{d-1} + d-1$. 
It is well-known that $\mu(\m^r) = {r + d-1 \choose d-1}$. It is easy
to verify 
by induction on $d$ that ${r + d-1 \choose d-1} > r^{d-1} + d-1$ for
all $d\geq 3$ and for all $r>1$.  Hence $r=1$.

We now prove the sufficiency. With the assumptions in the theorem it
is easy to see that the equation of minimal multiplicity holds 
 for all $d \geq 1$. 
Let $
\j = (t_{1}^{-1}, x_{1}^r t_{2} + t_{2}^{-1}, x_d t_{1}, 
     x_i t_{1} + x_{i+1}^r t_{2}; 1 \leq i \leq d-1)
$ (put $r=1$ for $d \geq 3$). 
Then $\j \n = \n^2$. In view of Remark~\ref{sally1},
$\b(\m,\m^{r})_{\n}$ is CM with minimal multiplicity. 
\qed

\bt
\label{two03}
Let $(R, \m)$ be a CM local ring of dimension $d \geq 2$.  Let 
$I$ be an $\m$-primary parameter ideal in $R$.  Then ${\cal
B}(\m,I)_{\n}$ is CM with minimal
multiplicity if and only if $R$ is a regular local ring and  $\ell
(I + \m^2/ \m^2) \geq d-1$.  
\et
\pf Suppose $\b(\m,I)_{\n}$ is CM with minimal multiplicity. Then 
 $e_{d-1}(\m|I) = \mu(I) - d + 1 = 1$ (Theorem~\ref{last}, 
Theorem~\ref{last1}).  This implies that $e(\m) \leq e_{d-1}(\m|I)  =
1$ and hence  $e(\m) =1$.   By a result of
Rees \cite{rees5},  there exists $x_{1}, \ldots, 
x_{d-1} \in I$ and $x_{d} \in \m$ such that 
$e_{d-1}(\m| I) = e(x_{1}, \ldots, x_d)=1$. 
Hence $\m=(x_{1}, \ldots, x_d)$ and  $\ell (I+\m^2/\m^2) \geq d-1$.

Conversely, since $R$ is a regular local ring
\beqn
    \mu(\n) - \dim~\b(\m,I)_{\n} + 1
= 2 + \mu(\m) + \mu(I) - (d+2) + 1
= d+1.
\eeqn
Since  $\ell(I + \m^2/\m^2) \geq d-1$, there exists a regular system of
parameters $x_{1}, \ldots,  x_d$ in $R$ such that $I = (x_{1}, \ldots,
x_{d-1}, x_d^r)$. This implies 
$1 \leq e_{q}(\m|I) \leq e(x_{1}, \ldots,  x_d) = 1$ for $q = 0,
\ldots d-1$ (Lemma~\ref{isw}). Thus
\beqn
     e(\n) 
=   \f{1}{2^d} \left[ 1 
    + \sum_{q=0}^{d-1} 2^{q+1}  
    + \sum_{q = 0}^{d-1} 2^{d-1-q} (d-q) \right]
=    d+1.
\eeqn
Hence $e(\n) = \mu(\n) - \dim~\b(\m,I)_{\n} + 1$. 
Let
$$
\j = (t_1^{-1}, x_1 t_2 + t_2^{-1}, x_d t_1, 
     x_i t_1 + x_{i+1} t_2; 1 \leq i \leq d-2,
     x_{d-1} t_1 + x_{d}^{r} t_2).
$$
Then $\j \n = \n^2$. In view of Remark~\ref{sally1}, $\b(\m,I)_{\n}$
 is CM with minimal multiplicity. 
\qed

\brm
{\em
In \cite[pg. 3015]{verma2}, J.K. Verma  asked the following
question:  
If $(R, \m)$
is a CM local ring, $I$ is any ideal in $R$, $\b(I)_{\n}$ is CM with
minimal multiplicity, then is it true that $R$ has minimal
multiplicity. This question does not have an affirmative answer in
general. The following example shows  that  there exist 
extended Rees algebras which are  CM with minimal multiplicity even
though $R$ does not  have minimal multiplicity. For details on this
example the author is requested to see 
\cite[Example~4.2.8, Example~4.2.9]{clare3}.
}
\erm

\bex
\label{ex2}
{\em 
Let $R= k[[x^4,x^5,x^7]]$ where $x$ is an indeterminate,
$\m=(x^4, x^5, x^7)$,  $I_{1}=(x^4,\m^2)$, $I_{2}=(\m^2)$. Then $R$ is
a CM ring which does not have minimal multiplicity, but
${\b}(I_{1},I_{2})_{\n}$   and ${\b}(I_{1})_{\n}$ are CM with minimal
multiplicity. \qed }  
\eex

\small{School of Mathematics, 
         Tata Institute of Fundamental Research,
         Homi Bhabha Road,}\\
         \small{Mumbai-400 005, India.}

{Chennai Mathematical Institute, 92, G. N. Chetty
Road, Chennai 600 017, India. \\
email: clare@cmi.ac.in}
\end{document}